\documentclass[reqno,11pt]{amsart}
\usepackage{latexsym,amssymb,amsfonts, amsmath}
\usepackage{a4wide}
\theoremstyle{definition}

\newcommand{\ee}{\mathbb{E}}

\newcommand{\rr}{\mathbb{R}}

\newcommand{\pp}{\mathbb{P}}

\newcommand{\ds}{\displaystyle}

\begin{document}
\bibliographystyle{xxx}

\mbox{}

\author[]{Fran\c cois Bolley}

\title[]{Optimal coupling for mean field limits}

\thanks{Ceremade, Umr Cnrs 7534, Universit\'e Paris-Dauphine, Place du Mar\'echal de Lattre de Tassigny,  F-75775 Paris cedex 16. E-mail address: bolley@ceremade.dauphine.fr}



\begin{abstract}
We review recent quantitative results on the approximation of mean field diffusion equations  by large systems of interacting particles, obtained by optimal coupling methods. These results concern a larger range of models, more precise senses of convergence and links with the long time behaviour of the systems to be considered.
\end{abstract}

\maketitle

Let us consider a Borel probability distribution $f_t = f_t(X)$ on $\rr^d$ evolving according to the McKean-Vlasov equation
\begin{equation}\label{mkv}
\frac{\partial f_t}{\partial t}  = \sum_{i,j=1}^d a_{ij} \frac{\partial^2 f_t}{\partial X_i \partial X_j} + \sum_{i=1}^d  \frac{\partial}{\partial X_i} \big( b_i[X, f_t] \, f_t \big),\qquad t>0, \; X \in \rr^{d}.
\end{equation}

Here $a = (a_{ij})_{1 \leq i, j \leq d}$ is a nonnegative symmetric $d \times d$ matrix; moreover, given $X$ in $\rr^d$ and a Borel probability measure $p$ on $\rr^d$,
$$
b_i[X,p] = \int_{\rr^d} b_i(X,Y) \, dp(Y), \qquad 1 \leq i \leq d
$$
where $b(X,Y) = (b_i(X,Y))_{1 \leq i \leq d}$ is a vector of $\rr^d.$ Equation~\eqref{mkv} has the
following natural probabilistic interpretation: if $f_0$ is a distribution on $\rr^d$, the solution $f_t$ of \eqref{mkv} is the law at time $t$ of the $\rr^{d}$-valued process $(X_t)_{t\geq 0}$ 
evolving according to the mean field stochastic differential equation
\begin{equation}\label{eds}
dX_t    = \sigma \, dB_t - b[X_t, f_t] \, dt.   
\end{equation}
Here the $d \times d$ matrix $\sigma$ satisfies $\sigma \sigma^* = 2 a$, $(B_t)_{t \geq 0}$ is a Brownian motion in $\rr^{d}$ and $f_t$ is the law of $X_t$ in $\rr^{d}$. It is a mean field equation in the sense that the evolution of $X_t$ is obtained by averaging the contributions $b(X_t, Y)$ over the system, according to the distribution $df_t(Y).$ Existence and uniqueness of solutions to~\eqref{mkv} and~\eqref{eds} are proven in~\cite{Mel96} for globally Lipschitz drifts $b$ and initial data $f_0$ with finite second moment. Non globally Lipschitz drifts are discussed in section~\ref{nonlip}.

\bigskip

Two instances of such evolutions are particularly interesting. First of all, when $\rr^d$ is the phase space of positions $x \in \rr^{d'}$ and velocities $v \in \rr^{d'}$ with $d = 2 d'$, one is interested in the Vlasov-Fokker-Planck equation
\begin{equation}\label{vfp}
\frac{\partial f_t}{\partial t} + v \cdot \nabla_x f_t - \big( \nabla_x U *_x \rho[f_t] \big)  \cdot \nabla_v f_t = \Delta_v f_t + \nabla_v \cdot ((A(v) + B(x)) f_t), \qquad t>0, \; x, v \in \rr^{d'}.
\end{equation}
Here $a \cdot b$ denotes the scalar product of two vectors $a$ and $b$ in $\rr^{d'}$, whereas $\nabla_v$, $\nabla_v \cdot$ and $\Delta_v$ respectively stand for the gradient, divergence and Laplace operators with respect to the velocity variable $v \in \rr^{d'}$. Moreover  $ \ds \rho[f_t](x)= \ds \int_{\rr^{d'}} f_t(x,v) \, dv $ is the macroscopic density in the space of positions $x \in \rr^{d'}$, or
the space marginal of $f_t$; $U=U(x)$ is an interaction potential in the position space and  $*_x$ stands for the convolution with respect to the position variable $x \in \rr^{d'};$ finally $A(v)$ and $B(x)$ are respectively friction and position confinement terms. This equation is used in the modelling of diffusive stellar matter (see~\cite{Dol99} for instance). 
\smallskip

We are also concerned with the space homogeneous equation
\begin{equation}\label{edpMKV}
\frac{\partial f_t}{\partial t}  = \Delta_v f_t + \nabla_v \cdot (( \nabla_v V + \nabla_v W *_v f_t)f_t), \qquad t>0, \; v \in \rr^{d}
\end{equation}
with $d=d'.$  Here $V$ and $W$ are respectively exterior and interaction potentials in the velocity space, and this equation is used in the modelling of space homogeneous granular media (see~\cite{BCCP98}). 
\bigskip

The {\it particle approximation} of~\eqref{mkv} consists in introducing $N$ processes $(X^{i,N}_t)_{t \geq 0}$, with $1 \leq i \leq N$, which evolve no more according to the drift  $b[X_t, f_t]$ generated by the distribution $f_t$ as in~\eqref{eds}, but according to its discrete counterpart, namely the empirical measure
\[
\hat{\mu}^N_t = \frac{1}{N} \sum_{i=1}^{N} \delta_{X^{i,N}_t}
\]
of the particle system $(X^{1,N}_t, \dots, X^{N,N}_t).$ In other words we let the processes $(X^{i,N}_t)_{t \geq 0}$ solve
\begin{equation}\label{edsXi}
dX_t^{i,N} = \sigma \, dB^i_t - \frac{1}{N} \sum_{i=1}^N b(X_t^{i,N}, X^{j,N}_t) \, dt, 
\qquad 1 \leq i \leq N.
\end{equation}
Here the $(B^i_t)_{t \geq 0}$'s are $N$ independent standard Brownian motions on $\rr^d$ and we assume that the initial data $X^{i,N}_0$ for $1 \leq i \leq N$ are independent variables with given law $f_0.$

\medskip

The mean field force $b[X_t, f_t]$ in~\eqref{eds} is replaced in~\eqref{edsXi} by the pairwise actions $\ds \frac{1}{N} b(X^{i,N}_t,X^{j,N}_t)$ of particle $j$ on particle $i.$ In particular, even in this case when the initial data $X^{i,N}_0$ are independent, the particles get correlated at all $t>0.$ But, since this interaction is of order $1/N$, it may be reasonable that two of these interacting particles (or a fixed number $k$ of them) become less and less correlated as $N$ gets large. 

In order to state this {\it propagation of chaos} property we let, for each $i \geq 1$,  $(\bar{X}^i_t)_{t \geq 0}$ be the solution of 
\begin{equation}
\label{eq:sdemkv}
\begin{cases}
  d\bar{X}^i_t    = \sigma \, dB^i_t - b[\bar{X}^i_t, f_t] \, dt\\
  \bar{X}_0^i=X_0^{i,N}
\end{cases}
\end{equation}
where $f_t$ is the distribution of $\bar{X}^i_t.$ The processes ${(\bar{X}^i_t)_{t \geq 0}}$ with $i\geq 1$ are independent since the initial conditions and driving Brownian motions are independent. Moreover they are identically distributed and their common law at time $t$ evolves according to~\eqref{mkv}, so is the solution $f_t$ of~\eqref{mkv} with initial datum $f_0.$ In this notation, and as $N$ gets large, we expect the $N$ processes $(X^{i,N}_t)_{t\geq 0}$ to look more and more like the $N$ independent processes ${(\bar{X}^i_t)_{t \geq 0}}$:

\medskip

\noindent
{\bf Theorem (\cite{Mel96}, \cite{Szn91})}  If $b$ is a Lipschitz map on $\rr^{2d}$ and $f_0$ a Borel distribution on $\rr^d$ with finite second moment, then, in the above notation, for all $t \geq 0$ there exists a constant $C$ such that
$$
\ee \vert X^{i,N}_t -  \bar{X}^i_t \vert^2 \leq \frac{C}{N}
$$
 for all $N.$
 
 This results also holds in a more general setting when the diffusion matrix can depend on $X$ and on the distribution in a Lipschitz way, and can be stated at the level of the paths of the processes on finite time intervals.
 
 First of all, it ensures that the common law $f^{1,N}_t$ of any (by exchangeability) of the particles $X^{i,N}_t$ converges
to $f_t$ as $N$ goes to infinity, according to
\begin{equation}\label{estloi}
W_2^2(f^{1,N}_t, f_t) \leq \ee \vert X^{i,N}_t -  \bar{X}^i_t \vert^2 \leq \frac{C}{N} \cdot
\end{equation}
Here the Wasserstein distance of order $p \geq 1$ between two Borel probability measures $\mu$ and $\nu$ on $\rr^q$ with finite moment of order $p \geq 1$ is defined by
$$
W_p (\mu, \nu) = \inf \big( \ee \vert X - Y \vert^p \big)^{1/p}
$$
where the infimum runs over all couples of random variables $(X,Y)$ with $X$ having law $\mu$ and $Y$ having law $\nu$ (see~\cite{Vil09} for instance).

Moreover, it proves a quantitative version of propagation of chaos : for all fixed $k$, the law $f^{k,N}_t$ of any (by exchangeability) $k$ particles $X^{i,N}_t$ converges
to the product tensor $(f_t)^{\otimes k}$ as $N$ goes to infinity, according to
$$
W_2^2(f^{k,N}_t, (f_t)^{\otimes k}) \leq \ee \vert (X^{1,N}_t, \dots, X^{k,N}_t) -  (\bar{X}^1_t, \dots, \bar{X}^k_t) \vert^2 \leq \frac{k C}{N} \cdot
$$

It finally gives the following first result on the convergence of the empirical measure $\hat{\mu}^N_t$ of the particle system to the distribution $f_t$: if $\varphi$ is a Lipschitz map on $\rr^d$, then
\begin{equation}\label{mesemp1}
\ee \Big\vert \frac{1}{N} \sum_{i=1}^N \varphi(X^{i,N}_t) - \! \int_{\rr^d} \! \! \varphi \, df_t \Big\vert^2 \leq 2 \, \ee \big\vert \varphi(X^{i,N}_t) - \varphi(\bar{X}^i_t) \big\vert^2 + 2 \, \ee \Big\vert  \frac{1}{N} \sum_{i=1}^N \varphi(\bar{X}^{i}_t) - \! \int_{\rr^d} \! \!  \varphi \, df_t \Big\vert^2 \! \leq \frac{C}{N}
 \end{equation}
by the Theorem and a law of large numbers argument on the independent variables $\bar{X}^{i}_t.$

\bigskip

Recent attention has been brought to improve these classical results in three directions:

1. enlarging the setting to non Lipschitz drift terms;

2. providing more precise estimates on the approximation of the solution to~\eqref{mkv} by the empirical measure of the particle system;

3. providing time uniform estimates when possible, in connexion with the long time behaviour of the solutions.

\section{Non Lipschitz drifts}\label{nonlip}

The interest for such mean field limits with locally but non globally Lipschitz drifts has been renewed by kinetic models, in which the interaction may become larger for larger relative velocities. 

For instance, for space homogeneous models such as~\eqref{edpMKV} with $W(z)=\vert z \vert^3$ and $d=1$, the difficulty brought by the non Lipschitz drift has been solved by convexity arguments, first in dimension one in~\cite{BRTV98}, then more generally in any dimension in~\cite{CGM06},~\cite{Mal01}. 

Space inhomogeneous collective behaviour models are under study in~\cite{BCC10}, for which known convexity arguments do not apply. 

\section{Deviation bounds for the empirical measure}

The averaged estimate~\eqref{mesemp1} ensures that the particle system is an appropriate approximation to solutions to~\eqref{mkv}.  However, when the particle system is used for numerical simulations, one may wish to establish estimates making sure that the numerical method has a very small probability to give wrong results.

This was achieved in~\cite{bgm09},~\cite{CGM06},~\cite{Mal01}, ~\cite{Mal03} by the use of (Talagrand) transportation inequalities. It is proved in these works, under diverse hypotheses on the initial data and in diverse contexts, that for all $t$ the law $f^{N,N}_t$ of the particle system at time $t$ satisfies a transportation inequality
$$
W_1(\nu, f^{N,N}_t)^2 \leq \frac{1}{c} H(\nu \vert f^{N,N}_t)
$$
for all measures $\nu$ on $\rr^{dN}$. Here $c$ may depend on $t$ (but neither on $N$ nor on $\nu$), and 
$$
H(\nu \vert f^{N,N}_t) = \int_{\rr^{dN}} \frac{d\nu}{d f^{N,N}_t} \, \ln \frac{d\nu}{df^{N,N}_t} \, df^{N,N}_t
$$
is the relative entropy of $\nu$ with respect to $f^{N,N}_t,$ to be interpreted as $+\infty$ if $\nu$ is not absolutely continuous with respect to $f^{N,N}_t.$ Then an argument by S.~Bobkov and F.~G\"otze, based on the Kantorovich-Rubinstein dual formulation 
\begin{equation}\label{kr}
W_1(\nu, f^{N,N}_t) = \sup \Big\{ \int_{\rr^{dN}} \Phi \, d\nu -  \int_{\rr^{dN}} \Phi \, df^{N,N}_t , \, \Phi \; 1{\textrm{-Lipschitz on}} \; \rr^{dN} \Big\}
\end{equation}
and the dual formulation of the entropy, ensures that 
$$
 \pp \Big[ \frac{1}{N} \sum_{i=1}^N \varphi(X^{i, N}_t) - \! \int_{\rr^d} \! \! \! \varphi \, df^{1,N}_t > r  \Big] \leq e^{- c \, N \, r^2}
$$
for all $N \geq 1, r >0$ and all $1$-Lipschitz maps $\varphi$ on $\rr^d$ (see~\cite[Chapter 22]{Vil09} for instance). Since moreover
\begin{eqnarray*}
\Bigl | 
\frac{1}{N} \sum_{i=1}^{N} \varphi (X_t^{i,N}) - \int_{\rr^d} \varphi\,df_t \Bigr | 
&\leq& 
\Bigl | 
\frac{1}{N} \sum_{i=1}^{N} \varphi (X_t^{i,N}) - \int_{\rr^d} \varphi\,df^{1,N}_t \Bigr | + W_1(f^{1,N}_t, f_t) \\
&\leq& 
\Bigl \vert \frac{1}{N} \sum_{i=1}^{N} \varphi (X_t^{i,N}) - \int_{\rr^d} \varphi\,df^{1,N}_t \Bigr | + \sqrt{\frac{C}{N}}
\end{eqnarray*}
by~\eqref{kr},~\eqref{estloi} and the bound $W_1 \leq W_2$, this ensures one-observable error bounds like
\begin{equation}\label{mesemp2}
\pp \left[ \Bigl | 
\frac{1}{N} \sum_{i=1}^{N} \varphi (X_t^{i,N}) - \int_{\rr^d} \varphi\,df_t \Bigr | > \sqrt{\frac{C}{N}} + r \right]
\leq 2\, e^{-c N r^2}
\end{equation}
for all $N, r$ and all $1$-Lipschitz maps $\varphi$ on $\rr^d$. In this argument we see how well adapted to this issue are the Wasserstein distances: from their definition they can easily be bounded from above by simple estimates on the processes, as in~\eqref{estloi}, and in turn they lead to straightforward estimates on Lipschitz observables, by~\eqref{kr}.

\medskip

Uniformly on Lipschitz observables, and for the Wasserstein distance $W_1$ which, up to moment conditions, metrizes the narrow topology on measures, estimates like
\begin{equation}\label{mesemp3}
\pp \Big[ W_1(\hat{\mu}^N_t, f_t) > r \Big] \leq e^{-\lambda N r^2},
\end{equation}
or even
$$
\pp \Big[ \sup_{0 \leq t \leq T} W_1(\hat{\mu}^N_t, f_t) > r \Big] \leq e^{-c N r^2}
$$
were reached in~\cite{BGV05}, provided $N$ is larger than an explicit $N_0(r),$ hence ensuring that the probability of observing any significant deviation during a whole time period $[0,T]$ is small. Also bounds were obtained on the {\it pointwise} deviation of a mollified empirical measure around the solution $f_t.$

\section{Time uniform estimates and long time behaviour}

Under convexity assumptions on the potentials $V$ and $W,$ the solution to the space homogeneous granular media equation~\eqref{edpMKV} has been proven to converge algebraically or exponentially fast to a unique steady state (see~\cite{CMV03} for an entropy dissipation proof based on interpreting~\eqref{edpMKV} as a gradient flow in the Wasserstein space, and also~\cite{CGM06},~\cite{Mal01},~\cite{Mal03}). In this setting one can hope for time uniform constants in estimates~\eqref{mesemp1}-\eqref{mesemp2}-\eqref{mesemp3}, which were obtained in~\cite{BGV05},~\cite{CGM06},~\cite{Mal01},~\cite{Mal03}. 

Also convergence to equilibrium (through a contraction argument in $W_2$ distance) and time uniform deviation bounds were obtained in~\cite{bgm09} for solutions to the (now space inhomogeneous) Vlasov-Fokker-Planck equation~\eqref{vfp}.


\medskip

\begin{thebibliography}{10}

\bibitem{BRTV98}
{\sc S. Benachour, B. Roynette, D. Talay and P. Vallois.}
\newblock Nonlinear self-stabilizing processes. {I}. {E}xistence, invariant
  probability, propagation of chaos.
\newblock {\em Stoch. Proc. Appl. 75}, 2 (1998), 173--201.

\bibitem{BCCP98}
{\sc D. Benedetto, E. Caglioti, J.~A. Carrillo and M. Pulvirenti.}
\newblock A non-{M}axwellian steady distribution for one-dimensional granular
  media.
\newblock {\em J. Statist. Phys. 91}, 5-6 (1998), 979--990.

\bibitem{BCC10}
{\sc F. Bolley, J.~A.~Ca\~nizo and J.~A.~Carrillo.}
\newblock Work in preparation (2010).

\bibitem{bgm09}
{\sc F. Bolley, A. Guillin and F. Malrieu.}
\newblock Trend to equilibrium and particle approximation for a weakly
 selfconsistent Vlasov-Fokker-Planck equation.
 \newblock Accepted for publication in {\em M2AN} (2010).

\bibitem{BGV05}
{\sc F. Bolley, A. Guillin and C. Villani.}
\newblock Quantitative concentration inequalities for empirical measures on non-compact spaces.
\newblock {\em Prob. Theor. Rel. Fields 137}, 3-4 (2007), 541--593.

\bibitem{CMV03}
{\sc J.~A.~Carrillo, R.~J.~McCann and C.~Villani.}
\newblock Kinetic equilibration rates for granular media and related equations:
  entropy dissipation and mass transportation estimates.
\newblock {\em Rev. Mat. Iberoamericana 19}, 3 (2003), 971--1018.

\bibitem{CGM06}
{\sc P. Cattiaux, A. Guillin and F. Malrieu.}
\newblock Probabilistic approach for granular media equations in the non uniformly case.
\newblock  {\em Prob. Theor. Rel. Fields 140},  1-2 (2008), 19--40.

\bibitem{Dol99}
{\sc J. Dolbeault.}
\newblock Free energy and solutions of the Vlasov-Poisson-Fokker-Planck system: external potential and confinement (large time behavior and steady states).  
\newblock {\em J. Math. Pures Appl. 9}, 78, 2  (1999), 121--157.

\bibitem{Mal01}
{\sc F. Malrieu.}
\newblock Logarithmic {S}obolev inequalities for some nonlinear {PDE}'s.
\newblock {\em Stoch. Proc. Appl. 95}, 1 (2001), 109--132.

\bibitem{Mal03}
{\sc F. Malrieu.}
\newblock Convergence to equilibrium for granular media equations and their
  {E}uler schemes.
\newblock {\em Ann. Appl. Probab. 13}, 2 (2003), 540--560.

\bibitem{Mel96}
{\sc S. M\'el\'eard.}
\newblock {Asymptotic behaviour of some interacting particle systems; McKean-Vlasov and Boltzmann models.}  
\newblock {Lecture Notes in Math. 1627, Springer, Berlin, 1996.}

\bibitem{Szn91}
{\sc A.-S. Sznitman.}
\newblock Topics in propagation of chaos.
\newblock {Lecture Notes in Math. 1464, Springer, Berlin, 1991.}
  
\bibitem{Vil09}
{\sc C. Villani.}
\newblock {\em Optimal transport, old and new}.
\newblock Grundlehren der math. Wiss. 338, Springer, Berlin, 2009.

\end{thebibliography}
\end{document}